\documentclass[12pt,a4paper]{amsart}
\usepackage[utf8]{inputenc}
\usepackage[english]{babel}

\usepackage[a4paper,inner=2.8cm,outer=2.5cm,top=2.8cm,bottom=2.5cm]{geometry}

\usepackage{amsmath,amssymb,amscd,amsthm,amsfonts}
\usepackage{mathtools}
\usepackage{mathrsfs}
\usepackage{stmaryrd}
\usepackage{tensor}
\usepackage{cases}
\usepackage{xfrac}
\usepackage{empheq}

\usepackage{booktabs}
\usepackage{tabularx}
\usepackage{makecell}
\usepackage{ragged2e}
\usepackage{array}
\usepackage{microtype}
\usepackage{enumerate}
\usepackage{mdwlist}
\usepackage{xspace}
\usepackage{comment}
\usepackage{lipsum}
\usepackage{fix-cm}
\usepackage{layout}
\usepackage{lpic}
\usepackage{anyfontsize}
\usepackage{color}
\usepackage{xcolor}

\usepackage{graphicx}
\usepackage{float}
\usepackage{pdfpages}
\usepackage{tikz}
\usepackage{tikz-cd}
\usetikzlibrary{calc,decorations.markings,positioning}
\usepackage[all]{xy}

\usepackage[bottom]{footmisc}
\usepackage[colorlinks=true,linkcolor=blue,citecolor=purple,filecolor=magenta,urlcolor=cyan]{hyperref}

\theoremstyle{plain}
\newtheorem{theorem}{Theorem}[section]

\newtheorem{corollary}[theorem]{Corollary}

\theoremstyle{definition}
\newtheorem{definition}[theorem]{Definition}

\newtheorem{example}[theorem]{Example}

\theoremstyle{remark}
\newtheorem{remark}[theorem]{Remark}

\allowdisplaybreaks

\usepackage[font=small,labelfont=bf]{caption}
\begin{document}
\title[Positivity in electrical networks, Grassmannians and symplectic groups]{Positivity in electrical networks, Grassmannians and symplectic groups}

\author[V. Gorbounov]{V.~Gorbounov}
\address{V.~G.: Faculty of Mathematics, National Research University Higher School of Economics, Usacheva 6, 119048 Moscow, Russia}
\email{vgorb10@gmail.com }

\author[A. Kazakov]{A.~Kazakov}
\address{A.~K.: 
Centre of Integrable Systems, P. G. Demidov Yaroslavl State University, Sovetskaya 14, 150003,
Yaroslavl; Faculty of Mathematics, National Research University Higher School of Economics; Steklov Mathematical Institute of Russian Academy of Sciences, Moscow, Russia}
\email{anton.kazakov.4@mail.ru}

\begin{abstract}
We study the relation between electrical positivity introduced by Lam and Pylyavskyy
\cite{L2} and total positivity of symplectic matrices arising from circular planar electrical networks. A particularly convenient framework for our purpose is provided by a  construction of an   embedding  of electrical network  into the positive Lagrangian Grassmannian, see Section \ref{el-to-gras}.
Using this technique, we systematically compare these two types of positivity and provide criteria for the positivity of decomposition parameters, see Section \ref{mian-sec_2}.  
In Section \ref{el-pos-factor} we also derive explicit rational formulas recovering the parameters of the decomposition from our symplectic matrix entries and minors, and prove uniqueness under explicit nonvanishing conditions.
\medskip


\medskip

\end{abstract}

\medskip

\maketitle

\tableofcontents

\hspace{0.1in}

\section{Electrical networks}
\begin{definition} 
An \textbf{electrical network} $\mathcal{E}(\Gamma, \omega)$ is a connected weighted planar graph $\Gamma(V, E)$ embedded into a disk and  satisfying the following conditions:
\begin{itemize}
    \item[\textbullet] The nodes are partitioned into inner nodes $V_I$ and boundary nodes $V_B$, with each boundary node lying on the boundary circle;
    \item[\textbullet] The boundary nodes are enumerated clockwise from $1$ to $|V_B| = n$; the inner nodes are enumerated arbitrarily, starting from $n+1$;
    \item[\textbullet] Every edge is assigned a positive weight $t_{ij}$, representing its conductivity.
\end{itemize}

\end{definition}

Each electrical network gives rise to the following natural problem. Suppose voltages $\textbf{U} = \{U_1, \dots, U_n\}$ are applied at the boundary nodes. These boundary voltages  admit a unique harmonic extension to the inner nodes, and induce  the corresponding boundary currents $\textbf{I} = \{I_1, \dots, I_n\}$, which  are determined via Kirchhoff's law:

\begin{align} \label{kir-equat}
    \begin{cases}
        \sum \limits_{j \in V}t_{ij}\bigr(U_i-U_j\bigl)=0, \ \forall i \in V_I \\
          I_k:=\sum \limits_{j \in V}t_{ij}\bigr(U_i-U_j\bigl), \ \forall i\in V_B
    \end{cases}
\end{align}
\begin{theorem} \textup{\cite{CM}} \label{resp_matr_sta}
    Let $\mathcal{E}(\Gamma, \omega)$ be an electrical network with $n$ boundary nodes. Then there exists  a  matrix $M_\mathcal{E}=(x_{ij}) \in \text{Mat}_{n \times n}(\mathbb{R})$ that relates $\textbf{U}$ and $\textbf{I}$ to each other  as follows:
$$M_\mathcal{E}\textbf{U}=\textbf{I}.$$
This matrix is called the \textbf{response matrix} of the network $\mathcal{E}$ and satisfies the following conditions:
\begin{itemize}
    \item[\textbullet] The matrix $M_\mathcal{E}$ is symmetric;
    \item[\textbullet] The sum of the entries in each row $($column$)$ is $0.$
\end{itemize}
\end{theorem}

\begin{figure}[H]
    \centering
    \includegraphics[scale=1.0]{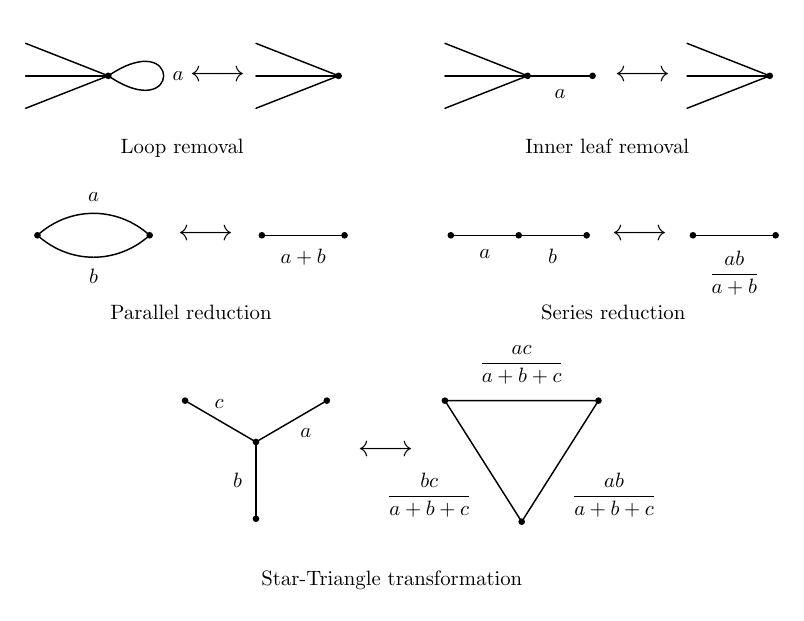}
\caption{Electrical transformations}
    \label{fig:el_trans}
\end{figure}

\begin{theorem} \textup{\cite[Theorem $1$]{CM}} \label{gen-th-1}
   Two electrical networks are equivalent under the transformations shown in Fig. \ref{fig:el_trans} if and only if they share the same response matrix. 
\end{theorem}

\begin{figure}[H]
     \hspace*{-7mm}
     \includegraphics[scale=1.0]{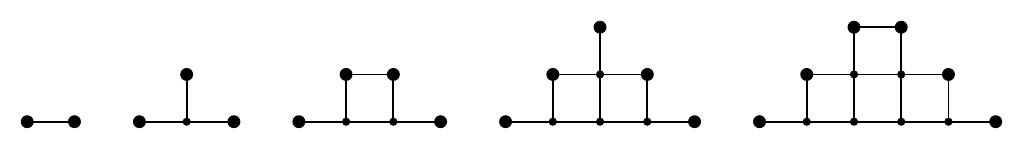}
     \caption{Standard networks $\mathcal{S}_n$, boundary nodes are bold } \label{fig:standart}
     \end{figure}

By Theorem \ref{gen-th-1}, we will consider electrical networks up to electrical transformations. We denote by $E_n$ the quotient set of electrical networks with $n$ boundary nodes. It is natural to consider generic electrical networks, i.e. networks that cannot be simplified by electrical transformations and have the maximum possible number of edges:

\begin{definition}   
\label{def:well-con} 
    An electrical network $\mathcal{E}(\Gamma, \omega) \in E_n$  is called \textbf{well-connected}  
    if it is equivalent to the standard network $\mathcal{S}_n$, see Fig. \ref{fig:standart}.
\end{definition}
It is also natural to extend this consideration to cactus electrical networks, which feature extended conductivity functions \(\omega \colon E(G)\rightarrow \mathbb{R}_{\geqslant 0}\cup \{+\infty \}\). According to Kirchhoff's laws \eqref{kir-equat}, an edge with zero conductivity can be deleted, while an edge with infinite conductivity can be contracted.  We denote by \(\overline{E}_{n}\) the quotient set of cactus electrical networks with \(n\) boundary nodes modulo electrical transformations.

\subsection{From electrical networks to non-negative Grassmannians} \label{el-to-gras}

Electrical network theory can naturally be viewed as  part of the theory of non-negative Grassmannians \cite{L}, \cite{Pos}.
\begin{definition}
The totally non-negative Grassmannian \(\mathrm{Gr}_{\geq 0}(k, m)\) is the subset of the real Grassmannian \(\mathrm{Gr}(k, m)\) consisting of all points whose Plücker coordinates \(\Delta _{I}\) (for \(I \in \binom{[m]}{k}\)) can be chosen to be non-negative. Similarly, the totally positive Grassmannian \(\mathrm{Gr}_{> 0}(k, m)\) is the subset consisting of all points whose Plücker coordinates \(\Delta _{I}\) can be chosen to be strictly positive.
\end{definition}

\begin{theorem} \textup{\cite[Theorem $3.3$]{BGKT}} \label{th: main_gr}
 Consider an  electrical network  $\mathcal{E}(\Gamma, \omega) \in E_n$ with a response matrix $M_{\mathcal{E}}=(x_{ij}).$ Then there is an embedding $\mathcal{I}: E_n\to \mathrm{Gr}_{\geq 0}(n-1,2n)$ such that $\mathcal{I}(\mathcal{E})$ is the row space of the $(n\times 2n)$ matrix
\begin{equation*}
\Omega_n(\mathcal{E})=\left(
\begin{array}{cccccccc}
x_{11} & 1 & -x_{12} & 0 & x_{13} & 0 & \cdots & (-1)^n\\
-x_{21} & 1 & x_{22} & 1 & -x_{23} & 0 & \cdots & 0 \\
x_{31} & 0 & -x_{32} & 1 & x_{33} & 1 & \cdots & 0 \\
\vdots & \vdots &  \vdots &   \vdots &  \vdots & \vdots & \ddots & \vdots 
\end{array}
\right)    
\end{equation*}
In particular, we have:
\begin{itemize}
    \item 
    The dimension of the row space of $\Omega_n(\mathcal{E})$ is equal to $n-1$;
    \item Every $n-1 \times n-1$ minor of  $\Omega_n(\mathcal{E})$  is non-negative.
\end{itemize}
\end{theorem}
\begin{theorem} \textup{\cite{BGKT}, \cite{L}}
   The embedding \(\mathcal{I}\) extends to the embedding \(\overline{\mathcal{I}}: \overline{E}_n \to \mathrm{Gr}_{\geq 0}(n-1, 2n)\), which can be described in purely combinatorial terms $($see \cite[Section 5.2]{L}$)$.
\end{theorem}

Let us define a \((2n - 2)\)-dimensional subspace \(V \subset \mathbb{R}^{2n}\) by the following conditions:
\[
V := \left\{v \in \mathbb{R}^{2n} \mid \sum_{i=1}^{n} (-1)^i v_{2i} = 0, \sum_{i=1}^{n} (-1)^i v_{2i-1} = 0\right\}.
\]

In the subspace \(V\), consider the basis:
\begin{gather*}
w_1 = (1, 0, 1, 0, \dots, 0, 0, 0), \; w_2 = (0, 1, 0, 1, \dots, 0, 0, 0),  \ \dots, \  w_{2n-2} = (0, 0, 0, 0, \dots, 1, 0, 1). 
\end{gather*}
One can  verify that the rows of the matrix \(\Omega_n(\mathcal{E})\) belong to the space \(V\). By choosing a specific basis for \(V\), we obtain an embedding of \(E_{n}\) into \(\mathrm{Gr}_{\geq 0}(n-1,2n-2)\). Certain bases of \(V\) are distinguished by their algebraic properties. Now we provide a special embedding \(\mathcal{L}: E_n \to \mathrm{Gr}_{\geq 0}(n-1, 2n-2)\) (which, of course, can be naturally extended to the embedding \(\overline{\mathcal{I}}: \overline{E}_n \to \mathrm{Gr}_{\geq 0}(n-1, 2n-2)\)).
\begin{theorem} \textup{\cite[Section $4.3$]{BKT}, \cite{Tal}}
    Let us suppose that $n$ is odd, then there exists a special basis $F$ of the subspace $V$ such that expanding the rows of the matrix $\Omega_n(\mathcal{E})$ with respect to this basis gives us a matrix $\overline{\Omega}_n(\mathcal{E})$ that defines a point of $\mathrm{LG}_{\ge 0}(n-1, 2n-2)$ with respect to the non-degenerate symplectic form:
\begin{equation} \label{form}
    \overline{\Lambda}_{2n-2} = \begin{pmatrix} \Lambda_n & 0 \\ 0 & -\Lambda_n \end{pmatrix},
\end{equation}
where $\Lambda_n \in \mathrm{Mat}_{n-1 \times n-1}$ is defined as follows:
\begin{equation} \label{form2}
\Lambda_n = \begin{pmatrix}
0 & 1 & 0 & 0 & \dots \\
-1 & 0 & -1 & 0 & \dots \\
0 & 1 & 0 & 1 & \dots \\
\vdots & \vdots & \vdots & \ddots & \vdots
\end{pmatrix}.
\end{equation}
In other words, we have an embedding $\mathcal{L}: E_n \to \mathrm{LG}_{\geq 0}(n-1, 2n).$
\end{theorem} 

\begin{figure}[H]
\hspace*{-7mm}
     \includegraphics[scale=1.0]{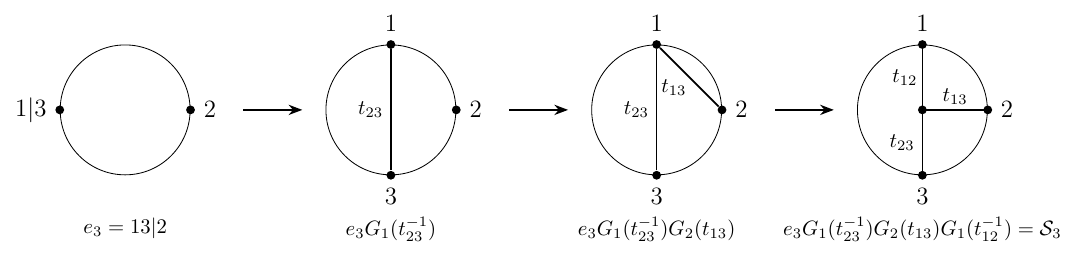}
     \caption{Adding spikes and bridges to $e_3=13|2$} \label{ex_3}
     \end{figure}

It is clear that if \(\Delta_{12\dots n-1}\Bigl(\overline{\Omega}_n(\mathcal{E})\Bigr) \neq 0\), then \(\overline{\Omega}_n(\mathcal{E})\) determines the same point in the Grassmannian as the following matrix:
\begin{equation}  \label{form_matrix}
   \bigr(\overline{\mathrm{E}}, A(\mathcal{E})\bigl) = \begin{pmatrix}
0 & \dots & 0 & 0 & 1 & a_{11} & a_{12} & a_{13} & \dots \\
0 & \dots & 0 & -1 & 0 & a_{21} & a_{22} & a_{23} & \dots \\
0 & \dots & 1 & 0 & 0 & a_{31} & a_{32} & a_{33} & \dots \\
\vdots & \ddots & \vdots & \vdots & \vdots & \vdots & \vdots & \ddots & \vdots
\end{pmatrix},
\end{equation}
 \begin{equation*}
\ \ \  \overline{\mathrm{E}}=\begin{pmatrix}
    0 & \dots & 0 & 0 & 1  \\
0 & \dots & 0 & -1 & 0  \\
0 & \dots & 1 & 0 & 0  \\
\vdots & \ddots & \vdots & \vdots & \vdots
\end{pmatrix} \in \mathrm{Mat}_{n-1 \times n-1}.
\end{equation*}
In turn, the matrix \(A(\mathcal{E})\) can be described in purely combinatorial electrical terms. 
\begin{definition}\label{positivity}
     A non-degenerate matrix \(X \in \operatorname{Mat}_{n \times n}(\mathbb{R})\) is called totally positive (non-negative) if all of its minors are strictly  positive (non-negative).
\end{definition}
\begin{theorem} \textup{\cite[Section $4.3$]{BKT}, \cite{Tal}} \label{well-con-th}
  Consider  a standard network $\mathcal{S}_n$. Then the matrix $A(\mathcal{S}_n) = (a_{ij}) \in \mathrm{Mat}_{n-1 \times n-1}(\mathbb{R})$ defined below is a totally non-negative symplectic matrix with respect to the skew-symmetric non-degenerate form \eqref{form2}:
\[
 \quad A(\mathcal{S}_n)\Lambda_n A(\mathcal{S}_n)^T = \Lambda_n.
\]
Moreover, the matrix $A(\mathcal{S}_n)$ can be decomposed into the product of the elementary totally non-negative \textbf{symplectic generators} $u_{i}(t) = E + t(E_{i+1,i} + E_{i-1,i})$:
\begin{equation} \label{prod_form}
    A(\mathcal{S}_{n}) = \prod_{1 \le i < j \le n} u_{j-i}(t_{ij}^{(-1)^{i+j}}),
\end{equation}
where $E$ is the identity matrix, $E_{ij}$ is the matrix unit, and $E_{n,n-1} = E_{0,1} = 0$; the product is taken over all pairs $(i,j)$ ordered in reverse lexicographical order; and
here $t_{ij}$ are  conductivities of  $\mathcal{S}_n$ indexed  according to the rule described in  \cite{Tal}.
\end{theorem}
\begin{proof}
    We provide an independent, combinatorial proof of this statement. Let us consider a hollow cactus network \(e_n := 1n \vert{} 2(n-1) \vert{} \dots\ \vert{} \frac{n-1}{2}+1 \), which is constructed as follows: consider a network with edges \(i(n-i+1)\) and set their conductivities to infinity, i.e., contract them (for instance, $e_3=13|2$ and \(e_5 = 15 \vert{} 24 \vert{} 3\)). 
   Now, consider the operators \(G_i(t) \in \mathrm{Mat}_{2n \times 2n}, i=1, \dots, n-1\) corresponding to adding a bridge for even indices and adding a spike for odd indices.  Matrices of these operators (see \cite[Section \(6\)]{BGKT}), denoted by the same symbols, have the following form:
    \begin{gather*}
G_1(t) = E + t(E_{1,2} - E_{1,2n}), \\ G_k(t) = E + t(E_{k,k-1} + E_{k-1,k}), \ k \in \{2, \dots, n-1\}.
\end{gather*}
It can be verified that the standard electrical network \(\mathcal{S}_{n}\) with conductivities \(t_{ij}\) (indexed according to the convention in \cite{Tal}) can be obtained by adding spikes and bridges as follows (see Fig. \ref{ex_3}):
\begin{equation} \label{orbit}
    \mathcal{S}_n=e_n\prod_{1 \le i < j \le n} G_{j-i}(t_{ij}^{(-1)^{i+j}}).
\end{equation}
Let \(F\) be a basis for the subspace \(V\) whose  vectors are the rows of the following matrix:
\begin{equation} \label{basis-main}  
 B_n:=\left( \begin{array} {ccccccccccccc}
0 & 0 & 0 & 0 &  \dots & \overset{n\mathrm{-th}}{1} & 0 & 1 & 0 & \dots & 0 &0 & 0 \\
0 & 0 & 0 &  0 & \dots & 0 & 1 & 0 & 1 & \dots & 0 & 0 & 0 \\
0 & 0 & 0 &  0 & \dots & 0 & 0 & 1 & 0 & \dots & 0 & 0 & 0 \\
\vdots & \vdots & \vdots &  \vdots & \ddots & \vdots & \vdots & \vdots & \vdots & \ddots & \vdots & \vdots & \vdots \\
0 & 0 & 0 &  0 & \dots & 0 & 0 & 0 & 0 & \dots & 1 & 0 & 1 \\
0 & 1 & 0 & 0 & \dots & 0 & 0 & 0 & 0 & \dots & 0 & 0 & -1 \\
1 & 0 & 1 & 0 & \dots & 0 & 0 & 0 & 0 &  \dots & 0 & 0 & 0 \\
0 & 1 & 0 & 1 & \dots & 0 & 0 & 0 & 0 &    \dots & 0 & 0 & 0 \\
0 & 0 & 1 & 0 & \dots &  0 & 0 & 0 & 0 &   \dots & 0 & 0 & 0 \\
\vdots & \vdots & \vdots &  \vdots & \ddots & \vdots & \vdots & \vdots & \vdots & \ddots & \vdots & \vdots & \vdots 
\end{array} \right).
\end{equation} 

One can verify that \(U_i(t)\) are the matrices of the operators \(G_i(t)\) in the defined basis \(F\):
\begin{gather*} \label{form_u}
 B_nG_i(t)=U_i(t)B_n, \   U_i(t):= \begin{pmatrix}
E & 0   \\
0 & u_i(t) 
\end{pmatrix}.
\end{gather*}
Additionally, the following holds:
$$\bigr(\overline{E}, E\bigl)B_n=\overline{\mathcal{I}}(e_n).$$
Indeed, using \cite[Section 5.2]{L}, we obtain that
$\Delta_{J}\bigl(\overline{\mathcal{I}}(e_n)\bigr)=1$
for each $J$ obtained by picking exactly one index from each non-degenerate bracket (i.e., excluding the brackets $(\tfrac{n-1}{2}+1, \tfrac{n-1}{2}+1)$ and $(2n, 2n)$)
\[
(1, 2n-1),\ (3, 2n-3),\ (5, 2n-5),\ \dots,\ (2, 2n-2),\ (4, 2n-4),\ (6, 2n-6), \ \dots;
\]
and $\Delta_{J}\bigl(\overline{\mathcal{I}}(e_n)\bigr)=0$ otherwise. The same identities hold for the Plücker coordinates of \(\bigl(\overline{E}, E\bigl)B_n\).

Combining all these results together, we obtain
\begin{align*}
\mathcal{L}(\mathcal{S}_n)=\mathcal{L}\Bigl(e_n\prod_{1 \le i < j \le n} G_{j-i}(t_{ij}^{(-1)^{i+j}})\Bigr)=\overline{\mathcal{L}}(e_n)\prod_{1 \le i < j \le n} U_{j-i}(t_{ij}^{(-1)^{i+j}})= \\=(\overline{E}, E)\prod_{1 \le i < j \le n} U_{j-i}(t_{ij}^{(-1)^{i+j}})=\Bigl(\overline{E}, \prod_{1 \le i < j \le n} u_{j-i}(t_{ij}^{(-1)^{i+j}})\Bigr)=\Bigl(\overline{E}, A(\mathcal{S}_{n})\Bigr). 
\end{align*}
To finish the proof, it remains to note that
$$U_i(t)\overline{\Lambda}_{2n-2}U^T_i(t)=\overline{\Lambda}_{2n-2}, \ (\overline{E}, E)\overline{\Lambda}_{2n-2}(\overline{E}, E)^T=0 \implies \Bigl(\overline{E}, A(\mathcal{S}_{n})\Bigr) \in \mathrm{LG}_{\geq 0}(n-1, 2n-2);$$
$$u_i(t)\Lambda_nu^T_i(t)=\Lambda_n \implies A(\mathcal{S}_{n})\Lambda_nA^T(\mathcal{S}_{n})=\Lambda_n.$$
\end{proof}
\begin{remark}
  From the proof of Theorem \ref{well-con-th}, it follows that for odd \(n\), the image of \(\overline{\mathcal{I}}\) is the closure of the orbit of \(e_{n}\) under the action of \(G_i(t), i \in \{1, \dots, n-1\}\) (compare with \cite[Theorem 5.11]{BGKT}).
  
Furthermore, the proof of Theorem \ref{well-con-th} provides a purely combinatorial interpretation of the vertex model described in \cite{GoT}, and precisely relates the embedding \(\overline{\mathcal{I}}\) to the semigroup \((El_{2n})_{\geq 0}\) defined in \cite[Proposition 4.1]{L2}.
\end{remark}

Finally, note that Theorem \ref{well-con-th} can be refined:
\begin{theorem} \label{main-cror}
Consider a well-connected network $\mathcal{E}$ with an odd number of boundary nodes. Then the matrix $\overline{\Omega}_n(\mathcal{E})$ defines a point in $\mathrm{LG}_{> 0}(n-1,2n-2)$.     
\end{theorem}
\begin{proof}
    To prove this statement, it is sufficient to demonstrate that the matrix $$A(\mathcal{E})=u_1(*)\bigl(u_2(*)u_1(*)\bigr)\bigl(u_3(*)u_2(*)u_1(*)\bigr)\dots$$ is totally positive. Let us consider the following correspondence between elementary generators and transpositions:
    $$x_i(t):=E+tE_{i, i+1} \to (i, i+1) \in S_{n-1} \ \ 
        y_i(t):=E+tE_{i+1, i} \to \overline{(i, i+1)} \in S_{n-1}.$$
   
   One can verify that, with respect to the defined correspondence, the product $$\bigl(u_2(*)u_1(*)\bigr)\bigl(u_3(*)u_2(*)u_1(*)\bigr)\dots$$ contains the highest permutations $$[n, n-1, \dots, 1]=(1, 2)\bigl((2, 3)(1, 2)\bigr)\bigl((3, 4)(2, 3)(1, 2)\bigr)\dots$$ as well as the same permutations with overlined indices. Thus, according to \cite[Factorization schemes]{Fomin}, \(A(\mathcal{E})\) is  totally positive. 
\end{proof}

\begin{remark}
Unlike the embedding \(\mathcal{I}\), the embedding \(\mathcal{L}\)  does not map
the cells of the electrical-network CW complex to the corresponding
positroid cells of the positive Grassmannian. 

Indeed, consider the case \(n=5\) and let \(\mathcal{E}\) be the network
obtained from the standard network by setting
\(t_{45}=+\infty\).
This operation corresponds to edge contraction, and, hence, leads to a lower-dimensional cell in the electrical-network CW complex. However, the permutation computation from the proof of Theorem \ref{main-cror}    demonstrates that
\(\mathcal{L}(\mathcal{E})
\in \mathrm{LG}_{>0}(4,8)\).
Thus, although \(\mathcal{E}\) belongs to a lower-dimensional cell of the
electrical-network CW complex, its image under \(\mathcal{L}\)
still lies in the strictly positive part of the Lagrangian Grassmannian.
\end{remark}

\section{Positivity in symplectic groups} \label{mian-sec_2}
The embedding $\mathcal{L}$ provides a natural framework for studying positivity in the symplectic group. A symplectic matrix \(A\) (with the respect to the form \eqref{form2}) is positive if and only if \((\overline{E}, A)\) determines a point of the non-negative Lagrangian Grassmannian, see \cite{Pos}. In turn, the properties of this point govern whether \(A\) admits a decomposition into elementary symplectic generators with positive parameters.
\begin{theorem} \label{main_prop}
A symplectic matrix \(A\) decomposes into the product \eqref{prod_form} with positive parameters if and only if the point \((\overline{E}, A)\) belongs to the image of the embedding $\mathcal{L}$.
\end{theorem}
 \begin{proof}
Consider a positive symplectic matrix \(A\) and a point \((\overline{E}, A) \in \mathrm{LG}_{\ge 0}(n-1, 2n-2)\). By Theorem \ref{well-con-th}, if this point lies in the image of the embedding \(\mathcal{L}\), then \(A\) decomposes into a product of elementary symplectic generators with positive parameters.

Conversely, suppose $A$ decomposes. Then the point of the form \eqref{form_matrix} belongs to the image of $\mathcal{L}$ by the physical interpretation of the   decomposition parameters from Theorem \ref{well-con-th}.    
 \end{proof}
 Notice that the factorizability of a matrix into the product \eqref{prod_form} with strictly positive parameters is equivalent to its factorizability into any other product -- parametrized by a reduced word of the longest element of \(S_{n}\) -- of elementary symplectic generators with strictly positive parameters. Indeed, from the proof of Theorem \ref{well-con-th} (see also Remark \ref{3-rem}), it follows that \(u_i(t)\) satisfies the relations from \cite[Theorem 4.1]{L2}; it remains to apply the standard arguments from \cite[Section 4.1]{L2}. Combining this observation with Theorem \ref{well-con-th} and Theorem \ref{main_prop}, we obtain the following results: 
\begin{corollary} \label{criteria}
  A symplectic matrix \(A\) admits a decomposition into a product -- parametrized by a reduced word for the longest element of \(S_{n}\) -- of elementary symplectic generators with  strictly positive parameter, if and only if $$(\overline{E}, A)B_n \in \mathrm{Gr}_{>0}(n-1, 2n).$$ 
\end{corollary}
\begin{remark}
   According to Corollary \ref{criteria}, verifying  the positivity of all decomposition parameters is equivalent to the positivity tests in \(\mathrm{Gr}(n-1, 2n)\) described in \cite{BGK}.
\end{remark}
From the considerations of continuity, we also obtain:
\begin{corollary} \label{criteria1}
  A symplectic matrix \(A\) admits a decomposition into a product  of elementary symplectic generators with  strictly positive parameter, if and only if $$(\overline{E}, A)B_n \in \mathrm{Gr}_{\geq 0}(n-1, 2n).$$ 
\end{corollary}
We now detail our technique using the following statement:
\begin{theorem} \label{contr}
The non-surjectivity of \(\mathcal{L}\) implies the existence of totally positive symplectic matrices that cannot be decomposed into elementary totally non-negative symplectic generators with positive parameters.
\end{theorem}
\begin{proof}
The embedding \(\mathcal{L}\) is already seen to be non-surjective when \(n=3\). Indeed, consider the point:
\begin{equation*}
Y:=    \begin{pmatrix}
0  & -1  & \frac{1}{2} & 1  \\
1 & 0  & 1 & 4 
\end{pmatrix}  \in \mathrm{LG}_{> 0}(2, 4).
\end{equation*}
Since not all diagonal elements of the right-hand block are greater than \(1\), this matrix cannot be decomposed into a product of elementary generators. Therefore, by Theorem \ref{main_prop}, the point \(Y\) does not belong to the image of \(\mathcal{L}\). 

Now we focus on the case n=5. Specializing \eqref{basis-main}, we obtain
\begin{equation*} \label{main-5-id}
\Omega_5(\mathcal{E})=\overline{\Omega}_5(\mathcal{E})B_5, \ \ 
 B_5=\begin{pmatrix}
0 & 0 & 0 & 0 & 1 & 0 & 1 & 0 & 0 & 0 \\
0 & 0 & 0 & 0 & 0 & 1 & 0 & 1 & 0 & 0 \\
0 & 0 & 0 & 0 & 0 & 0 & 1 & 0 & 1 & 0 \\
0 & 0 & 0 & 0 & 0 & 0 & 0 & 1 & 0 & 1 \\
0 & 1 & 0 & 0 & 0 & 0 & 0 & 0 & 0 & -1 \\
1 & 0 & 1 & 0 & 0 & 0 & 0 & 0 & 0 & 0 \\
0 & 1 & 0 & 1 & 0 & 0 & 0 & 0 & 0 & 0 \\
0 & 0 & 1 & 0 & 1 & 0 & 0 & 0 & 0 & 0 
\end{pmatrix}.
\end{equation*} 
We show that  the following  totally positive symplectic matrix (see \cite{Link}) cannot be decomposed into any product of elementary  symplectic generators with positive parameters:
\begin{equation} \label{contr_example}
   A=\begin{pmatrix}
 \frac{10}{3} & \frac{80}{3} & \frac{8}{3} & \frac{7}{3} \\
 1 & 9 & 1 & 1 \\
 4 & 41 & 6 & 13 \\
\frac{1}{3} & \frac{11}{3} & \frac{2}{3} & \frac{7}{3}
\end{pmatrix} 
\end{equation}


Let us consider the point \(X \in \mathrm{LG}_{> 0}(4, 8)\) associated with \(A\):
\begin{equation*}
X:=\begin{pmatrix}
0 & 0 & 0 & 1 & \frac{10}{3} & \frac{80}{3} & \frac{8}{3} & \frac{7}{3} \\
0 & 0 & -1 & 0 & 1 & 9 & 1 & 1 \\
0 & 1 & 0 & 0 & 4 & 41 & 6 & 13 \\
-1 & 0 & 0 & 0 & \frac{1}{3} & \frac{11}{3} & \frac{2}{3} & \frac{7}{3}
\end{pmatrix}
=
\begin{pmatrix}
3 & 0 & -2 & 0 & 1 & 7 & 0 & -5 \\
-3 & 0 & 1 & 0 & 0 & 2 & 1 & 6 \\
2 & 1 & 1 & -1 & -1 & -2 & 1 & 5 \\
6 & 1 & 2 & 0 & 0 & 1 & 0 & -3
\end{pmatrix}. 
\end{equation*}

One can verify that $XB_5 \notin \mathrm{Gr}_{\geq 0}(4,10): \Delta_{1234}(XB_5)\Delta_{12310}(XB_5)<0$. Therefore the matrix \eqref{contr_example}
can not be decomposed into the product \eqref{prod_form} with strictly positive parameters.

Suppose that the matrix \(A\) can be decomposed into another product of elementary symplectic generators with strictly positive parameters, i.e., \(A=\prod u_i(t_i), \forall t_i >0\). Then the following identities hold:
\begin{gather*}
\bigr(\overline{E}, A\bigl)B_5=\Bigl(\bigr(\overline{E}, E\bigl)\prod U_{i}(t_i)\Bigl)B_5=XB_5 \notin \mathrm{Gr}_{\geq 0}(4,10), \ \   U_i(t_i)= \begin{pmatrix}
E & 0   \\
0 & u_i(t_i) 
\end{pmatrix}.  
\end{gather*}
On the other hand, let us consider the following matrices:
\begin{align*}
G_1(t) = E + t(E_{1,2} - E_{1,10}), \
G_2(t) = E + t(E_{2,1} + E_{2,3}), \\
G_3(t) = E + t(E_{3,2} + E_{3,4}), \
G_4(t) = E + t(E_{4,3} + E_{4,5}).
\end{align*}

The following identities  can be verify directly (see \cite{Link}):
\[
\forall i:  B_5 G_i(t) = U_i(t)B_5.
\]
At the same time, $$XB_5=\bigr(\overline{E}, E\bigl)B_5 G_{i}(t_i)\in \mathrm{Gr}_{\geq 0}(4, 10).$$ Indeed, \((\overline{E}, E)B_5\) belongs to \( \mathrm{Gr}_{\geq 0}(4, 10)\), and that each \(G_{i}\) preserves total non-negativity in the Grassmannian. The latter holds because  \(G_{i}\) coincide with the first four matrices of a bridge- or spike-adding operators from \cite[Proposition $5.2$]{BGKT}. This contradiction demonstrates that matrix $A$ cannot be decomposed into any product of elementary symplectic generators with strictly positive parameters.
\end{proof}
\begin{remark}
   For \(n=5\), the embedding \(\mathcal{L}\) can be written in terms of the response matrix entries as follows:
    \begin{equation*}
\begin{pmatrix}
-x_{14} - x_{15} & 0 & x_{15} & 0 & 1 & -x_{12} - x_{13} - x_{14} - x_{15} & 0 & x_{13} + x_{14} + x_{15} \\
x_{24} + x_{25} & 0 & -x_{25} & 0 & 0 & -x_{12} & 1 & -x_{23} - x_{24} - x_{25} \\
-x_{34} - x_{35} & 1 & x_{35} & -1 & -1 & x_{13} & 1 & -x_{13} - x_{23} \\
-x_{14} - x_{24} - x_{34} & 1 & -x_{45} & 0 & 0 & -x_{14} & 0 & x_{14} + x_{24}
\end{pmatrix}
    \end{equation*}
\end{remark}

\section{Electrical positivity} \label{el-pos-factor}
In this section we will make a comment on the positivity in the particular electrical group which appeared above. The electrical Lie algebras and electrical Lie groups were introduced in \cite{L2} and generalized in \cite{BG1}. The electrical positivity was investigated in a general situation in the upcoming paper \cite{BG2}. Here we will solve the inverse problem in dimension $4$. The solution to the inverse problem in arbitrary dimension is straightforward. 
This can be used to investigate the true meaning of this particular type of electrical positivity in the spirit of \cite{BG2}.
\begin{theorem}\label{inversion}
Let
\[
\mathfrak{u}_1=E_{21},\qquad
\mathfrak{u}_2=E_{12}+E_{32},\qquad
\mathfrak{u}_3=E_{23}+E_{43},\qquad
\mathfrak{u}_4=E_{34},
\]
and put as earlier $u_k(t)=\exp(tu_k)=I+tu_k$.
Suppose that $M=(m_{rs})\in\mathrm{SL}_4(\mathbb R)$
admits a factorization
\[
M=u_1(a)u_2(b)u_3(c)u_4(d)
  u_1(e)u_2(f)u_3(g)
  u_1(h)u_2(i)u_1(j).
\]
Note that this word is written in the reverse order compared to Section \ref{mian-sec_2}.

For increasing index sets $I,J$ of equal cardinality,
write $\Delta_{I,J}=\det M_{I,J}$. Define
\[
\begin{aligned}
C_{12}&=1+m_{24}-m_{44},&
C_{13}&=m_{14},\\
C_{14}&=m_{44}-1,&
C_{15}&=m_{34}-m_{14},\\
C_{23}&=\Delta_{34,34}-m_{13}-1,&
C_{24}&=\Delta_{12,34},\\
C_{25}&=m_{33}-m_{13}-1,&
C_{34}&=\Delta_{12,23}-\Delta_{12,34}-m_{42},\\
C_{35}&=m_{32}-m_{12},&
C_{45}&=1+m_{31}-m_{11}.
\end{aligned}
\]
If
\[
C_{13}C_{14}C_{15}C_{24}C_{25}C_{35}\ne0,
\]
then the decomposition in the displayed order is unique,
and its parameters are
\[
\begin{aligned}
a&=\frac{C_{12}}{C_{13}},&
b&=\frac{C_{13}}{C_{14}},&
c&=\frac{C_{14}}{C_{15}},&
d&=C_{15},\\
e&=\frac{C_{23}C_{14}}{C_{13}C_{24}},&
f&=\frac{C_{24}C_{15}}{C_{14}C_{25}},&
g&=\frac{C_{25}}{C_{15}},\\
h&=\frac{C_{34}C_{25}}{C_{24}C_{35}},&
i&=\frac{C_{35}}{C_{25}},&
j&=\frac{C_{45}}{C_{35}}.
\end{aligned}
\]
\end{theorem}
\begin{proof}
Put
\[
U_k(t)=I+tu_k,\qquad
(u_1,u_2,u_3,u_4)
=(E_{21},E_{12}+E_{32},E_{23}+E_{43},E_{34}),
\]
and suppose
\[
M=U_1(a)U_2(b)U_3(c)U_4(d)
  U_1(e)U_2(f)U_3(g)U_1(h)U_2(i)U_1(j).
\]
Write $\Delta_{I,J}=\det M_{I,J}$, and define
\[
\begin{aligned}
C_{12}&=1+m_{24}-m_{44},&
C_{13}&=m_{14},\\
C_{14}&=m_{44}-1,&
C_{15}&=m_{34}-m_{14},\\
C_{23}&=\Delta_{34,34}-m_{13}-1,&
C_{24}&=\Delta_{12,34},\\
C_{25}&=m_{33}-m_{13}-1,&
C_{34}&=\Delta_{12,23}-\Delta_{12,34}-m_{42},\\
C_{35}&=m_{32}-m_{12},&
C_{45}&=1+m_{31}-m_{11}.
\end{aligned}
\]

All factors following $U_4(d)$ fix the fourth coordinate
vector. Hence the fourth column gives
\[
m_{14}=bcd,\quad m_{24}=cd+abcd,\quad
m_{34}=d+bcd,\quad m_{44}=1+cd.
\]
Consequently,
\[
a=\frac{C_{12}}{C_{13}},\qquad
b=\frac{C_{13}}{C_{14}},\qquad
c=\frac{C_{14}}{C_{15}},\qquad d=C_{15}.
\]

Remove these four factors by setting
\[
N=u_4(-d)u_3(-c)u_2(-b)u_1(-a)M.
\]
Then
\[
N=u_1(e)u_2(f)u_3(g)u_1(h)u_2(i)u_1(j),
\]
so its third column satisfies
\[
N_{13}=fg,\qquad N_{43}=g,\qquad N_{23}-N_{43}=efg.
\]
The corresponding row operations give
\[
N_{13}=\frac{C_{24}}{C_{14}},\qquad
N_{43}=\frac{C_{25}}{C_{15}},\qquad
N_{23}-N_{43}=\frac{C_{23}}{C_{13}}.
\]
For example,
\[
N_{13}=(1+ab)m_{13}-bm_{23}
=\frac{m_{13}m_{24}-m_{14}m_{23}}{m_{44}-1}.
\]
Thus
\[
e=\frac{C_{23}C_{14}}{C_{13}C_{24}},\qquad
f=\frac{C_{24}C_{15}}{C_{14}C_{25}},\qquad
g=\frac{C_{25}}{C_{15}}.
\]

Finally, set
\[
P=u_3(-g)u_2(-f)u_1(-e)N
 =u_1(h)u_2(i)u_1(j).
\]
Multiplication and substitution yield
\[
P_{12}=i=\frac{C_{35}}{C_{25}},\qquad
P_{22}-1=hi=\frac{C_{34}}{C_{24}},\qquad
P_{31}=ij=\frac{C_{45}}{C_{25}}.
\]
Therefore
\[
h=\frac{C_{34}C_{25}}{C_{24}C_{35}},\qquad
i=\frac{C_{35}}{C_{25}},\qquad
j=\frac{C_{45}}{C_{35}}.
\]
These formulas hold wherever the displayed denominators
are non zero, as they are for the matrix under consideration.
\end{proof}
\begin{remark} \label{3-rem}
The formulas of Theorem \ref{inversion} determine the decomposition parameters for
every reduced word of the longest element \(w_0\in S_5\). Indeed, any two
such words are connected by commutations and braid moves. Factors with
\(\lvert i-j\rvert>1\) commute, while for \(\lvert i-j\rvert=1\) the
electrical star--triangle relation \cite{BG1}, \cite{GoT}, \cite[Section $4.1$]{L2} gives 
\[
  u_i(x)u_j(y)u_i(z)
  = u_j(yz/D)\,u_i(D)\,u_j(xy/D),
  \qquad D=x+z+xyz.
\]
This follows from the argument for an electrical group introduced in \cite{BG1}. Namely, inside a root copy $SL_2$ the electrical generator is a conjugate of the standard Chevalley nilpotent and the generators of the standard electrical group from \cite{BG1} coincide with these in Theorem \ref{inversion} except for the 
last one, which is a standard Chevalley nilpotent.

Thus each braid move replaces \((x,y,z)\) by \((yz/D,D,xy/D)\).
Composing these moves with the formulas of Theorem~3.1 expresses the
new parameters as rational functions of the same matrix entries and
minors. Each move is defined whenever \(D\ne0\) and preserves strict
positivity.
\end{remark}
\begin{example}
Consider the totally positive symplectic matrix $A$ defined by \eqref{contr_example}. 
We will recover the parameters in the decomposition using Theorem \ref{inversion}.

\[
\begin{aligned}
C_{12}&=-\frac13,&
C_{13}&=\frac73,&
C_{14}&=\frac43,&
C_{15}&=\frac{32}{3},&
C_{23}&=\frac53,\\
C_{24}&=\frac13,&
C_{25}&=\frac73,&
C_{34}&=-\frac43,&
C_{35}&=\frac{43}{3},&
C_{45}&=\frac53.
\end{aligned}
\]
Consequently,
\[
\begin{aligned}
a&=\frac{C_{12}}{C_{13}}=-\frac17,
&
b&=\frac{C_{13}}{C_{14}}=\frac74,
&
c&=\frac{C_{14}}{C_{15}}=\frac18,
&
d&=C_{15}=\frac{32}{3},\\
e&=\frac{C_{23}C_{14}}{C_{13}C_{24}}=\frac{20}{7},
&
f&=\frac{C_{24}C_{15}}{C_{14}C_{25}}=\frac87,
&
g&=\frac{C_{25}}{C_{15}}=\frac7{32},\\
h&=\frac{C_{34}C_{25}}{C_{24}C_{35}}=-\frac{28}{43},
&
i&=\frac{C_{35}}{C_{25}}=\frac{43}{7},
&
j&=\frac{C_{45}}{C_{35}}=\frac5{43}.
\end{aligned}
\]

Direct multiplication verifies that
{\small{
\[
A=
u_1\!\left(-\frac17\right)
u_2\!\left(\frac74\right)
u_3\!\left(\frac18\right)
u_4\!\left(\frac{32}{3}\right)
u_1\!\left(\frac{20}{7}\right)
u_2\!\left(\frac87\right)
u_3\!\left(\frac7{32}\right)
u_1\!\left(-\frac{28}{43}\right)
u_2\!\left(\frac{43}{7}\right)
u_1\!\left(\frac5{43}\right).
\]}}
All denominators in the recovery formulas are nonzero.
Hence, by Theorem \ref{inversion}, this decomposition is unique in the given
specified order. In particular, although every minor of \(A\)
is strictly positive, two of its decomposition parameters,
\(a\) and \(h\), are negative, which agrees with the earlier calculation.
\end{example}
\medskip
\noindent\textbf{Acknowledgments.}  We thank  Dmitry~Talalaev and Yevgeniy~Pavlov for valuable discussions. We express our gratitude to the anonymous reviewer whose detailed comments  helped significantly  improve the quality of the paper.
 Anton Kazakov gratefully acknowledges support from the Basic Research Program of HSE University (HSE-BR-2025-062). Vassily Gorbounov gratefully acknowledges support from the Basic Research Program of HSE University (HSE-BR-2025-023).

\end{document}